\documentclass[11pt]{article}
\usepackage{tikz}
\usetikzlibrary{shapes,calc, arrows, matrix, positioning, fit, calc}
\usepackage{amsmath,bm}
\usepackage{graphicx}
\usepackage{authblk}

\usepackage{placeins}
\usepackage{booktabs}
\usepackage[T1]{fontenc}
\usepackage{chemfig}
\usepackage{algorithm,algorithmic}
\usetikzlibrary{matrix}
\usepackage{float}
\usepackage{subfigure}
\usepackage{subcaption}
\graphicspath{ {Figures/} }

\usepackage{fullpage}  
\usepackage{hyperref} 
\hypersetup{colorlinks=true,linkcolor=blue,citecolor=blue}
\usepackage{cleveref} 
\usepackage{url} 
\usepackage{lscape}

\title{A high-order accurate moving mesh finite element method for the radial Kohn--Sham equation}

\author{Zheming Luo}
\author{Yang Kuang\thanks{Corresponding author. Email: {\tt ykuang@gdut.edu.cn}}}
\affil{School of Mathematics and Statistics \&  Center for Mathematics and Interdisciplinary Science (CMIS), Guangdong University of Technology,  China}
\date{\today}

\begin{document}
\maketitle

\begin{abstract}
In this paper, we introduce a highly accurate and efficient numerical solver for the radial Kohn--Sham equation. The equation is discretized using a high-order finite element method, with its performance further improved by incorporating a parameter-free moving mesh technique. This approach greatly reduces the number of elements required to achieve the desired precision. In practice, the mesh redistribution involves no more than three steps, ensuring the algorithm remains computationally efficient. Remarkably, with a maximum of $13$ elements, we successfully reproduce the NIST database results for elements with atomic numbers ranging from $1$ to $92$.
\end{abstract}


\section{Introduction}
Kohn--Sham density functional theory (KSDFT) is the most commonly used electronic structure method in condensed matter physics and a widely adopted approach in quantum chemistry \cite{kohn1965self, hohenberg1964inhomogeneous}. In KSDFT, the electronic structure of a system is obtained by solving the Kohn--Sham equations, which pose a nonlinear eigenvalue problem. Despite being a fundamental task, solving these equations presents numerical challenges. One challenge arises from the singularities in the external potential near the nucleus, which cause the wavefunctions to vary smoothly between atoms but exhibit sharp changes close to the nucleus \cite{fiolhais2008primer}. A common approach to addressing numerical challenges in the core region is the use of pseudopotentials \cite{pickett1989pseudopotential}, where solving the atomic Kohn--Sham equation is crucial. Beyond pseudopotential generation \cite{oliveira2008generating}, the atomic Kohn--Sham equation is also essential for computing atomic properties \cite{kotochigova1997local} and constructing methods such as linearized augmented plane waves (LAPW) \cite{singh2006planewaves} and linearized muffin-tin orbitals (LMTO) \cite{skriver2012lmto}.

Various numerical methods have been developed to solve the radial Kohn--Sham equation, including the finite difference method \cite{romanowski2007numerical}, finite element methods \cite{romanowski2008bspline, certik2013dftatom, lehtola2019fully, lehtola2020fully, certik2024highorder}, the B-spline method \cite{romanowski2011bspline}, and the integral equation approach \cite{uzulis2022radial}, etc. Of these, the high-order finite element method stands out for its exponential convergence with respect to polynomial order \cite{certik2024highorder} and its mesh flexibility. Currently, most finite element methods employ either a uniform mesh or a radial mesh. Uniform meshes are inefficient in the core region due to the high element density required, whereas radial meshes offer more flexibility. The most commonly used radial mesh is the exponential mesh, which, through careful adjustment of its parameters, can achieve high accuracy with relatively few elements. However, despite this flexibility, radial meshes can be cumbersome to work with due to the numerous user-defined parameters involved.

In this work, we propose a moving mesh finite element method for solving the radial Kohn--Sham equation, which achieves high accuracy with a minimal number of elements by utilizing a moving mesh strategy, without requiring user-defined parameters. The strategy is based on the equidistribution principle introduced by de Boor \cite{de1973good}, which involves placing mesh points so that a measure of the solution error is evenly distributed across each subinterval. This principle has proven highly effective for moving mesh methods \cite{tang2005moving, huang2011adaptive}. 

The key ingredient of the moving mesh method is the monitor function, which controls the redistribution of the mesh, typically defined by the singularity of the problem. In this work, the monitor function is designed as the sum of the curvatures of the wavefunctions for the radial Kohn--Sham equation. This leads to a moving mesh equation, which in its discretized form becomes a tridiagonal linear system, and can be solved efficiently using methods like the Thomas algorithm. The solution of this equation provides the new mesh distribution. Since the Kohn--Sham equation is nonlinear, a precise update of the solution is unnecessary—interpolating the solution from the old mesh to the new one is sufficient. The new mesh is then used to solve the Kohn--Sham equation, and the moving mesh process is repeated until the solution converges. In practice, starting from a uniform mesh, convergence can be achieved within three iterations.

To achieve high accuracy, the radial Kohn--Sham equation is discretized using a high-order finite element method. Numerical experiments on the iron atom demonstrate the superiority of this approach over lower-order methods. For example, using a uniform mesh with polynomial order $p = 3$, achieving a precision of $10^{-6}$ Hartree requires 4600 elements. In contrast, the same accuracy is achieved with only 119 elements when $p = 10$, highlighting the necessity of employing high-order methods.

To further enhance efficiency, the moving mesh strategy is applied. This technique adaptively concentrates mesh points where higher resolution is needed, significantly reducing the number of elements. For instance, with $p = 3$, the number of elements required to achieve $10^{-6}$ accuracy decreases from 4600 to 143 when using a moving mesh. Similarly, with $p = 10$, the element count drops from 119 to just 10, illustrating the effectiveness of the moving mesh approach.

The computational performance is further improved by solving the generalized eigenvalue problem using the locally optimal block preconditioned conjugate gradient (LOBPCG) method \cite{Adaptive-finite-element3}. To accelerate convergence, we adopt a preconditioner following the approach in \cite{bao2012hadaptive}. This preconditioner consists of a discretized Laplacian and a shift operator with respect to the mass matrix, which dramatically reduces the number of LOBPCG iterations from thousands to just a few, significantly improving the computational efficiency.

We validate our method by comparing the results with the NIST  database \cite{kotochigova2009atomic} for elements ranging from hydrogen (H) to uranium (U), covering atomic numbers $Z = 1$ to $Z = 92$. The numerical results show excellent agreement with the NIST data for both total and orbital energies. Remarkably, these results are obtained using no more than 13 elements, with polynomial order $p = 10$, demonstrating the accuracy and efficiency of the proposed method.

The rest of this paper is structured as follows. Section 2 introduces the radial Kohn--Sham equation and the finite element method. Section 3 outlines the moving mesh strategy. Section 4 presents the numerical results, and we conclude in Section 5.

\section{Kohn--Sham equation}
\subsection{ Radial Kohn--Sham equation}
According to the density functional theory, the ground state energy of the atomic system of charge $Z$ with $N$ occupied orbitals can be obtained by solving the following three-dimensional Kohn--Sham equation.
\begin{equation}
\label{eq:KS}
\begin{aligned}
\left\{
\begin{array}{ll}
    \left(-\frac{1}{2}\nabla^{2}+\textit{V}_{\mathrm{eff}}(\bm{r})\right)\psi_{i}(\bm{r})=\varepsilon_{i}\psi_{i}(\bm{r}), ~i=1,2,...,N, \\
    \\ \int_\Omega\psi_{i}(\bm{r})\psi_{j}(\bm{r})d\bm{r}=\delta_{ij}, ~i,j=1,2,...,N,
\end{array}
\right.
\end{aligned}
\end{equation}
The effective potential $V_{\mathrm{eff}}(\bm{r})$ in above equation is expressed as 
\begin{equation}
\label{eq:veff}
 V_{\mathrm{eff}}(\bm{r})=V_{\mathrm{ext}}(\bm{r})+V_{\mathrm{Har}}[\rho(\bm{r})]+V_{\mathrm{xc}}[\rho(\bm{r})],
\end{equation}
where  $\rho(\bm{r})$  denotes the electron density, which can be written as 
\begin{equation*}
\rho(\bm{r})=\sum_{i=1}^{N}\vert\psi_{i}(\bm{r})\vert^2.
\end{equation*}
The first term in the effective potential is the electrostatic potential caused by the nucleus, which can be written as
 \begin{equation*}
 \label{eq:vext}
 V_{\mathrm{ext}}(\bm{r})=-\frac{Z}{\vert \bm{r}\vert},
 \end{equation*}
The second term in \eqref{eq:veff} is the Hartree potential describing the interaction potential between electrons, which can be written as 
 \begin{equation*}
 V_{\mathrm{Har}}(\bm{r})=\int_{\Omega}\frac{\rho(\bm{r}^{'})}{\vert\bm{r}-\bm{r}^{'}\vert}d\bm{r}^{'}.
 \end{equation*}
It is noted that  the Hartree potential is also the solution of the Poisson equation
\begin{equation}
\label{eq:poisson}
\nabla^{2}V_{\mathrm{Har}}(\bm{r})=-4\pi \rho(\bm{r}).
\end{equation}
The last term in \eqref{eq:veff} represents the exchange-correlation potential, which incorporates the effects of the Pauli exclusion principle (exchange) and non-classical electron-electron correlations beyond the classical Coulomb interaction.

Note that for the atomic system, it is spherically symmetric, therefore the Kohn--Sham equation can be further simplified to the radial form. Hereafter, we denote $(r,\theta,\phi)$ as the spherical coordinate, where $r$ is the radial distance, $\theta$ is the ploar angle, and $\phi$ is the azimuthal angle. For a single atom, the nucleus assumes to be positioned at the origin, and therefore the effective potential becomes
\begin{equation*}
V_{\mathrm{eff}}(\bm{r})=V_{\mathrm{eff}}(\vert \bm{r}\vert)=V_{\mathrm{eff}}({r}).
\end{equation*}
The wave function can be written as the product of the radial function $R_{nl}(r)$ and the angular function $Y_{lm}(\theta,\phi)$
\begin{equation}
\label{eq:wave function}
\psi(\bm{r})=\psi_{nlm}(\bm{r})=R_{nl}(r)Y_{lm}(\theta,\phi),
\end{equation}
where $n$, $l$ and $m$ are the main quantum numbers, angular quantum numbers and magnetic quantum numbers, respectively. In virtue of the spherical symmetric property, the radial Kohn--Sham equation  for the atom with charge $Z$ is derived: 
\begin{equation}
\label{eq:radial-ks}
-\frac{1}{2}(rR_{nl}(r))^{''}+\left(V_\mathrm{eff}(r)+\frac{l(l+1)}{2r^2})\right)rR_{nl}(r)=\varepsilon_{nl}rR_{nl}(r),
\end{equation}
where $\varepsilon_{nl}$ is the eigenvalue, and the radial function $R_{nl}(r)$ satisfies 
$\int_{0}^{\infty}r^{2}R_{nl}^{2}(r)dr=1.$ 
By introducing the substitution $P_{nl}(r)=rR_{nl}(r)$, \Cref{eq:radial-ks} can be rewritten as 
\begin{equation}
    \label{eq:radial-ks2}
    -\frac{1}{2}P_{nl}^{''}(r) + \left(
    V_{\mathrm{eff}}(r)+\frac{l(l+1)}{2r^2}
    \right)P_{nl}(r) = \varepsilon_{nl}P_{nl}(r),
\end{equation}
and the normalization condition is
\begin{equation*}
    \int_0^{\infty} P_{nl}^2(r)\mathrm{d}r = 1.
\end{equation*}
With the radial wavefunction, the radial electron density can be written as  
\begin{equation}
\label{eq:updateRho}
    \rho(r) = \frac{1}{4\pi}\sum_{nl}f_{nl}\frac{P_{nl}^2(r)}{r^2},
\end{equation}
where $f_{nl}$ is  electronic occupation, and the electron density is normalized as 
\begin{equation}
\label{eq:normal-rho}
    \int_0^{\infty} 4\pi \rho(r)r^2\mathrm{d}r = Z.
\end{equation}

Now we rewrite each term in the effective potential in the spherical coordinate. The external potential takes the form $V_\mathrm{ext}(r)=-Z/r$. The Poisson equation \eqref{eq:poisson} for the  Hartree potential now becomes
\begin{equation}
\label{eq:VH}
  \frac{1}{r^{2}}(r^{2}V_{\mathrm{Har}}^{'})^{'}=V_{\mathrm{Har}}^{''}(r)+\frac{2}{r}V_{\mathrm{Har}}^{'}(r)=-4\pi \rho(r).
\end{equation}
As to the exchange-correlation potential, we adopt the local density approximation (LDA). Specifically, the slater exchange potential and the Vosko-Wilk-Nusair (VWN4) correlation potential \cite{vosko1980accurate} are adopted. The exchange-correlation potential energy can be written as
\begin{equation*}
\epsilon_{\mathrm{xc}}^{LD}(\rho)=\epsilon_{x}^{LD}(\rho)+\epsilon_{\mathrm{c}}^{LD}(\rho),
\end{equation*}
where $\epsilon_{x}^{LD}$ is the exchange energy and $\epsilon_{\mathrm{c}}^{LD}$ is correlation energy. Their detailed formulations are
\begin{align*}
&\epsilon_{\mathrm{x}}^{LD}(\rho)=-\frac{3}{4\pi}(3\pi^{2}n)^\frac{1}{3},\\
&\epsilon_{\mathrm{c}}^{LD}(\rho)=\frac{a}{2}\{\log(\frac{x^{2}}{Y(x)})+\frac{2b}{q}B(x) -\frac{bx_{o}}{Y(x_{0})}[\log(\frac{(x-x_{0})^2}{Y(x)})+\frac{2(b+2x_{0})}{q}B(x)]\},
\end{align*}
where the following coefficients and functions were used in the above equations
\begin{align*}
&a = 0.0621814, \quad b = 3.72744, \quad c = 12.8352, \\
&x = \sqrt{r_{s}}, \quad r_{s} = \left( \frac{3}{4\pi \rho} \right)^{\frac{1}{3}}, \quad q = \sqrt{4c - b^{2}}, \\
&Y(y) = y^{2} + by + c, \quad B(x) = \arctan \left( \frac{q}{2y + b} \right), \quad y_{0} = -0.10498.
\end{align*}

With the above expressions, we can now deliver the formula for the total ground state energy of the system. The total energy of radial Kohn--Sham equation is given by
\begin{equation*}
E_{\mathrm{tol}}=E_{\mathrm{k}}+E_{\mathrm{Har}}+E_{\mathrm{xc}}+E_{\mathrm{ext}}.
\end{equation*}
The detailed expressions are given below. The kinetic energy is expressed as
   \begin{equation*}
   E_\mathrm{k}=\sum_{nl}\int\psi^{*}(r)(-\frac{1}{2}\nabla^{2})\psi(r)dr=\sum_{nl}f_{nl}\varepsilon_{nl}-4\pi \int_{0}^{\infty} V_{\mathrm{eff}}(r)r^{2}\rho(r)dr.
   \end{equation*}
The Hartree potential energy  is
\begin{equation*}
\label{33}
E_{\mathrm{Har}}=\int2\pi V_{\mathrm{Har}}(r)\rho(r)r^{2}dr,
\end{equation*}
   the exchange correlation energy is
   \begin{equation*}
   E_{\mathrm{xc}}=\int 4\pi\epsilon_{\mathrm{xc}}\rho(r)r^{2}dr,
   \end{equation*}
 and  the external potential energy is
   \begin{equation*}
 E_{\mathrm{ext}}=4\pi \int -\frac{Z}{r}\rho(r)r^{2}dr=-4\pi Z\int \rho(r)rdr.
   \end{equation*}

\subsection{Finite element discretization}
We adopt the finite element method to discretize the radial Kohn--Sham equation \Cref{eq:radial-ks2}. Based on the exponential decay of wave function and electron density, we truncate the infinite domain $[0,\infty)$ to a bounded region $\Omega=[0, R ]$ as the computational domain.  Denote the Hamiltonian operator as $\hat{H}:=-1/2{d^2}/{dr^2}+V_\mathrm{eff}(r)$, then the variational form of the radial Kohn--Sham equation can be expressed as: Find
$(\varepsilon_{i},p_{i})\in R\times H_{0}^{1}(\Omega)$, $i=1,\dots,N$,  such that
\begin{equation*}
 \label{eq:KS weak form}
 \begin{aligned}
\left\{
\begin{aligned}
&\int_{\Omega} v(r)\hat{H}P(r)dr=\varepsilon  \int_{\Omega} v(r)P(r)dr,\quad \forall v(r)\in H_{0}^{1}(\Omega)\\
&\int_{\Omega}P_{i}(r)P_{j}(r)dr=\delta_{ij}, \quad i,j=1,2,...,N,
\end{aligned}
\right.
\end{aligned}
\end{equation*}
where $v(r)$ is the test function. The computational domain $\Omega$ is then partitioned into $N_\mathrm{ele}$ elements. Denote the element set as  $\tau=\{\tau_{K},K=1,2,3,...,n_\mathrm{ele}\}$, then the finite element basis $\{\phi_l\}, ~{l=1,\dots,n_\mathrm{bas}}$ can be formed. Assume the finite element space $V_h$ be established on this partition, then we get the discrete expression of the continuous variational form on $V_{h}$: Find $(\varepsilon_{i}^{h},p_{i}^{h})\in R\times V_h, i=1,2,..,N$ such that
 \begin{equation}
 \label{eq:FEM}
 \left\{
 \begin{aligned}
&\int_{\Omega} \phi(r) \hat{H}P^{h}(r)dr=\varepsilon\int_{\Omega} \phi(r) P^{h}(r)dr,\quad \forall \phi\in H_{0}^{1}(\Omega)\\
&\int_{\Omega}P_{i}^{h}(r)(r)P_{j}^{h}(r)dr=\delta_{ij},\quad i,j=1,2,...,N,
\end{aligned}\right.
\end{equation} 
where  $P^h(r)$ is finite element approximation of $P(r)$ on the space $V_{h}$ with the form
  \begin{equation*}
P^{h}(r)=\sum_{l=1}^{n_\mathrm{bas}}P_{l}^h\phi_{l}(r).
  \end{equation*}
In the above equation, $P_{l}^h$ also represents the value of the radial function on the $l$-th degree of freedom. Due to the arbitrariness of $\phi$, we can substitute it with the basis function $\phi_{l}, l = 1, 2,\dots,n_\mathrm{bas}$ in \Cref{eq:FEM}. As a result, the following generalized eigenvalue problem is obtained
\begin{equation}
\label{eq:ks-gevp}
HX=\varepsilon MX,
\end{equation}
where  $H\in {R}^{n\times n}$ is the discrete Hamiltonian matrix, and $M$ is the mass matrix. The entries for $H$ and $M$ are 
\begin{equation}
\label{eq:hamil}
H_{i,j}=\frac{1}{2}\int_{\Omega}\nabla\phi_{j}\nabla\phi_{i}dr+\int_{\Omega}(\frac{l(l+1)}{2r^{2}}+V_\mathrm{eff}(r))\phi_{j}\phi_{i}dr,
\end{equation}
\begin{equation*}
M_{ij}=\int_{\Omega}\phi_{j}\phi_{i}dr.
\end{equation*}
Since the wavefunctions decay exponentially, the zero boundary condition can be applied for the discretized eigenvalue problem \eqref{eq:ks-gevp} when $\omega$ is large enough.

Similarly, the discretized linear system for the Hartree potential, as derived from \eqref{eq:VH}, can be expressed as:
\begin{equation}
    \label{eq:vh-ls}
    A V_\mathrm{Har} = b,
\end{equation}
where the entries of the matrix $A$ and the right-hand side vector $b$ are given by:
\begin{equation*}
    A_{i,j} = \int_\Omega \nabla \phi_j \cdot \nabla \phi_i - \frac{2}{r} \phi_{i}(r) \nabla \phi_j(r) \, dr,~
    b_{j} = \int_\Omega 4\pi \rho(r) \phi_{j}(r) \, dr.
\end{equation*}
Since the Hartree potential exhibits a $1/r$ decay behavior, a zero boundary condition cannot be directly applied. By substituting \eqref{eq:normal-rho} into \eqref{eq:VH}, we obtain:
\begin{equation*} 
    \lim_{r\to \infty} r^{2} \nabla V_{\mathrm{Har}}(r) = -Z,
\end{equation*}
which is equivalent to:
\begin{equation*}
    \nabla V_{\mathrm{Har}}(r) \sim -\frac{Z}{r^2}, \quad r \to \infty.
\end{equation*}
Integrating both sides yields:
\begin{equation*}
    V_{\mathrm{Har}}(r) \sim \frac{Z}{r}, \quad r \to \infty.
\end{equation*}
Thus, we derive the asymptotic property of $V_{\mathrm{Har}}(r)$. Consequently, we adopt $V_\mathrm{Har}(r) = {Z}/{r}$ as the boundary condition at the right end of the domain $\Omega$.

As $r$ approaches zero, we can analyze the asymptotic behavior of the Hartree potential by expanding the electron density as $\rho(r) = \sum_{i=0}^{\infty} c_{i} r^{i}$. Substituting this expansion into \eqref{eq:VH}, we find:
\begin{equation*}
    (r^{2} V_{\mathrm{Har}}'(r))' = -4\pi \sum_{i=1}^{\infty} c_{i} r^{i+2}.
\end{equation*}
Integrating this equation and requiring that $V_{\mathrm{Har}}(0)$ remains finite leads to:
\begin{equation}
    V_{\mathrm{Har}}'(r) = -4\pi \sum_{i=0}^{\infty} c_{i} \frac{r^{i+1}}{i+3}.
    \label{eq:VH_boundary}
\end{equation}
This indicates that the leading term behaves linearly, such that $V_{\mathrm{Har}}'(r) \propto r$ as $r \to 0$. Therefore, integrating \eqref{eq:VH_boundary} it follows
\begin{equation*}
    V_{\mathrm{Har}}(r) = -4\pi \sum_{i=0}^{\infty} c_{i} \frac{r^{i+2}}{(i+2)(i+3)} + C,
\end{equation*}
where the constant $C = V_{\mathrm{Har}}(0)$ is determined by Coulomb's law:
\begin{equation*} 
    V_{\mathrm{Har}}(0) = 4\pi \int_0^\infty r \rho(r) \, dr.
\end{equation*}
Finally, we find that $V'_{\mathrm{Har}}(0) = 0$, which will serve as the boundary condition at the left end of the computational domain $\Omega$.

\subsection{SCF iteration}
From the radial Kohn--Sham equation \eqref{eq:radial-ks2}, we can find that the Hamiltonian is determined by the electron density, while the electron density depends on the wavefunctions which are the solutions to radial KS equation. As a result, the radial KS equation is nonlinear. To resolve the non-linearity, the self-consistent field (SCF) iteration is adopted. Specifically, we use a linear mixing scheme 
\begin{equation*}
    \rho_\mathrm{new}(r) =  \alpha  \rho_\mathrm{in} + (1-\alpha)\rho_\mathrm{out},
\end{equation*}
where $\alpha$ is mixing parameter which is set as $0.618$ in our simulations. Notably, since \Cref{eq:radial-ks2} relies on the angular number $l$, it means that for an atom with different $L$ angular numbers, we have to solve  $L$ eigenvalue problems. The SCF method for the radial KS equation is outlines in \Cref{alg:scf}. 
\begin{algorithm}[H]
\caption{SCF method for radial Kohn--Sham equation.\label{alg:scf}}
\begin{algorithmic}[1]
    \STATE Input: computational domain  $\Omega$, partition $\tau$, atom configuration with electron number $Z$, initial electron density $\rho_{0}(r)$, tolerance $tol$
    \STATE $\rho_\mathrm{new} = \rho_0(r)$;
    \WHILE{$|E_\mathrm{new}-E_\mathrm{old}| > tol$}
        \STATE $E_\mathrm{old} = E_\mathrm{new}$, $\rho_\mathrm{in} = \rho_\mathrm{new}$;
        \STATE Construct effective potential $ V_{\text{eff}} $ from $ \rho_\mathrm{in}(r) $;
        \STATE Construct the generalized eigenvalue problems \eqref{eq:ks-gevp} for different $l$;  
        \STATE Solve eigenvalue problems to get $ \{\varepsilon_{nl},P_{nl}\} $;
        \STATE Update the electron density $\rho_\mathrm{out}$ from \eqref{eq:updateRho};
        \STATE Density mixing $ \rho_\mathrm{new}= \alpha \rho_\mathrm{out}(r) + (1-\alpha) \rho_{in}(r) $;
        \STATE Calculate the energy $E_\mathrm{new}$;
    \ENDWHILE
    \STATE Output $E_\mathrm{new}$.
\end{algorithmic}
\end{algorithm}

In order to reduce the number of self-consistent field iterations, we use the Thomas-Fermi approximation \cite{oulne2011variation} to generate the initial density and potential:
 \begin{align*}
 &V(r)=-\frac{Z_{\mathrm{eff}}(r)}{r},\\
 &Z_{\mathrm{eff}}(r)=Z(1+\alpha\sqrt{x}+\beta xe^{-\gamma\sqrt{x}})^{2}e^{-2\alpha \sqrt{x}},\\
 &x=r(\frac{128Z}{9\pi^{2}})^\frac{1}{3},\\
 &\alpha=0.7280642371,\quad \beta=-0.5430794693,\quad\gamma=0.3612163121,
 \end{align*}
and the corresponding initial charge density is given as
\begin{equation*}
\label{45}
\rho_\mathrm{0}(r)=-\frac{1}{3\pi^{2}}(-2V(r))^\frac{3}{2}.
\end{equation*}

\section{Moving mesh method}
To enhance the accuracy of the algorithm, it is crucial to distribute mesh points effectively on the domain.  Due to the rapid oscillations of the wavefunctions near the nucleus \cite{gulans2014exciting} and their smooth, slowly varying behavior at distances farther from the nucleus, a uniform mesh necessitates a significantly large number of points to maintain accuracy. Hence, a uniform mesh is not well-suited for solving the radial KS equation. The exponential mesh, which concentrates mesh points near the Coulomb singularity, is preferred and has been successfully implemented in advanced atomic structure codes such as \texttt{dftatom} \cite{certik2013dftatom}, \texttt{exciting} \cite{gulans2014exciting}, \texttt{HelFEM} \cite{lehtola2019fully}, and \texttt{featom} \cite{certik2024highorder}. While the exponential mesh effectively provides high-quality grids for the radial Kohn--Sham equation, it requires user-specified parameters, and the optimal parameters vary depending on the atom in question.  This parameter dependency limits the generality of these methods. By contrast, the moving mesh method eliminates the need for parameter tuning, offering a unified framework that applies uniformly to all atoms. 

The moving mesh method is based on the equidistribution principle, first introduced by de Boor \cite{de1973good} for solving boundary value problems in ordinary differential equations. This principle redistributes mesh points so that the solution error, measured in a certain way, is equalized across all subintervals. The error measure is typically described by a monitor function, which is able to capture the singularities of the problem and can be determined by geometric properties of the solution, such as arc length or curvature. In this work, we design the monitor function as
\begin{equation}
\label{eq:monitor}
M = \sqrt{\alpha + \sum_{nl} \left(\frac{\mathrm{d}P_{nl}}{\mathrm{d}x}\right)^2},
\end{equation}
where $\alpha$ is a positive parameter, set to 0.01 in this work. This monitor function ensures that mesh points are distributed such that the total arc length of all wavefunctions is equal in each subinterval. As a result, the mesh becomes denser near the origin, where the wavefunctions exhibit strong oscillations, and sparser in regions farther away, where the wavefunctions are smoother. This adaptive distribution is expected to enhance accuracy compared to a uniform mesh.

In the simulation, the moving mesh process is carried out with the help of an auxiliary interval $[0,1]$, which is uniformly partitioned and referred to as the computational domain \cite{tang2005moving}. The domain $[0,R]$, where the equation is defined, is known as the physical domain. Let $x, \xi$ denote the coordinates on the physical and  computational domains, respectively. The one-to-one coordinate transformation between these domains is given by 
\begin{equation*}
\label{eq:ctop}
    \left\{
        \begin{aligned}
            & x = x(\xi),~\xi \in [0,1],\\
            & x(0) = 0, ~ x(1) = R.
        \end{aligned}\right.
\end{equation*}
As a result, the equidistribution principle indicates the following form 
\begin{equation}
\label{eq:fordis}
    \int_{x(0)}^{x(\xi)} M \mathrm{d}x = \xi \int_0^RM \mathrm{d}x,
\end{equation}
holds true at each points $\xi_i = {i}/{(N+1)}, i=0,1,\cdots,N+1$ on the computational domain, where $N$ is the number of inner mesh points. After taking derivatives twice with respective to $\xi$, \Cref{eq:fordis} becomes
\begin{equation}
    \label{eq:eqdis}
    (Mx_\xi)_\xi = 0.
\end{equation}
The new mesh distribution can then be obtained by solving  \Cref{eq:eqdis}, along with the boundary conditions $x(0)=0$ and $x(1)=R$. 

Note that the monitor function $M$ depends on the wavefunctions, and these wavefunctions are functions of $x$. Consequently, \Cref{eq:eqdis} is a nonlinear equation, requiring a linearized method to solve it. We use the following linearization strategy
\begin{equation*}
    \label{eq:eqdis-linear}
    \left(M(x^{j})x_{\xi}^{j+1}\right)_\xi = 0, ~j=0,1,\cdots.
\end{equation*}
This leads to a semi-implicit scheme of the form:
\begin{equation*}
    M(x_{i+\frac{1}{2}}^j)(x_{i+1}^{j+1}-x_i^{j+1}
    -M(x_{i-\frac{1}{2}}^j)(x_{i}^{j+1}-x_{i-1}^{j+1}) = 0,
\end{equation*}
where $x_{i+1/2}=(x_i+x_{i+1})/2$. By combining this equation with the boundary conditions, a tridiagonal linear system is obtained, which can be solved efficiently by Thomas algorithm. 

From the construction of the monitor function and the solution of the associated equation \Cref{eq:eqdis}, it is evident that this process depends exclusively on the numerical solution of the wavefunctions $\psi_k$ computed on the current mesh. This observation leads to two key points: first, the process does not require parameter tuning, making it uniformly applicable to all atoms; second, due to the numerical errors inherent in the initial uniform mesh, a few iterations of the moving mesh process are necessary to achieve convergence.The algorithm is outlined in   \Cref{alg:mmm}. It is noted that in Step 7, the wavefunctions on the old mesh will interpolate to the new mesh and serve as the initial guess for the radial KS equation on the updated mesh.
\begin{algorithm}[H]
\caption{Moving mesh method for solving the radial KS equation \label{alg:mmm}}
\begin{algorithmic}[1]
    \STATE A uniform mesh $\{x_i^0\}$ on $[0,R]$, initial electron density, $tol$, $E_\mathrm{old}=0$. 
    
    \STATE Solve the radial KS equation to obtain $E_\mathrm{new}$, $\{\varepsilon_{nl}^\mathrm{new},P_{nl}^\mathrm{new}\}$.    
    \WHILE{$|E_\mathrm{new}-E_\mathrm{old}|>tol$}
        \STATE $E_\mathrm{old}=E_\mathrm{new}$.
        \STATE Generate the monitor function \Cref{eq:monitor}.
        \STATE Solve \Cref{eq:eqdis} to obtain the new mesh distribution $\{x_i^\mathrm{new}\}$.
        \STATE Interpolate wavefunctions from old mesh to new mesh.
        \STATE Solve the radial KS equation to obtain $E_\mathrm{new}$, $\{\varepsilon_{nl}^\mathrm{new},P_{nl}^\mathrm{new}\}$.
    \ENDWHILE
    \STATE Output $E_\mathrm{new}$.
\end{algorithmic}
\end{algorithm}

\section{Numerical Experiments}
In this section, we present numerical examples to illustrate the effectiveness of the proposed method. We begin by a systematic numerical study of the convergence rate and the impact of the moving mesh method on an iron atom, while also discussing the preconditioner used in the LOBPCG method. Then, we perform an experiment on a uranium atom to validate the method’s accuracy. Finally, we present and compare a set of numerical results with the NIST database \cite{kotochigova2009atomic}. In this context, atomic units are employed, specifically using Hartree for energy and Bohr for length.

In the following simulations, the SCF iteration and the moving mesh process ends up when the difference of the adjacent energies is below $10^{-8}$ Hartree.

\subsection{Numerical results of iron atom}
We first study the iron atom, whose nuclear charge is $Z=26$, and an electronic configuration of $1s^22s^22p^63s^23p^63d^64s^2$. As a result, three eigenvalue problems must be solved for the angular quantum numbers $l=0,1,2$. The orbitals required for each respective eigenvalue problem are $4, 2$, and $1$. The computational domain is set as $[0,20]$. In the following subsections, we will first demonstrate the effectiveness of the preconditioner in the LOBPCG method for solving the eigenvalue problem. Next, we will highlight the accuracy and efficiency of the moving mesh method.

\subsubsection{Preconditioner for the LOBPCG method}
To solve the generalized eigenvalue problem, we utilize the LOBPCG method \cite{Adaptive-finite-element3}, which has proven to be effective in electronic structure calculations \cite{castro2006octopus, gonze2009abinit}. One of the key strengths of the LOBPCG method is its use of a three-term recurrence relation, involving the current eigenvector approximation, the preconditioned residual, and the previous update direction. This allows for significantly reduced memory usage compared to Krylov subspace methods, making it suitable for large-scale problems. An important consideration for the LOBPCG method is the choice of preconditioner, which can greatly reduce the number of iterations and computation time. In this work, we construct the preconditioner following the approach in \cite{bao2012hadaptive}. 

It has the form of $T = L/2-\lambda M$, where $L/2$ is the discretized kinetic operator in Hamiltonian, and $\lambda$ is an approximated eigenvalue. For each eigenpair we will construct a preconditioner to accelerate the calculation, as a  result, the preconditioners are designed as 
\begin{equation} \label{eq:precon}
T_{l}^{(i)}=\left\{\begin{aligned}
\frac{1}{2} L-\lambda_l^{(i)} B & \text {, if } \lambda_{l}^{(i)} <0, \\
I & \text {, otherwise }
\end{aligned} \quad \text { for } l=1, \ldots, p.\right.
\end{equation}
In the practical simulations, the precondition process  involves the solution of linear systems $WT^{(i)}=(T_l^{(i)})^{-1}W^{(i)}$, which can be implemented by using the BiCG method, the same method employed to solve the Hartree potential.  Specifically, there is no need to accurately solve this linear system in the precondition process, and few iteration steps are performed in this process.

\begin{figure}[!ht]
    \centering
    \includegraphics[width=0.45\linewidth]{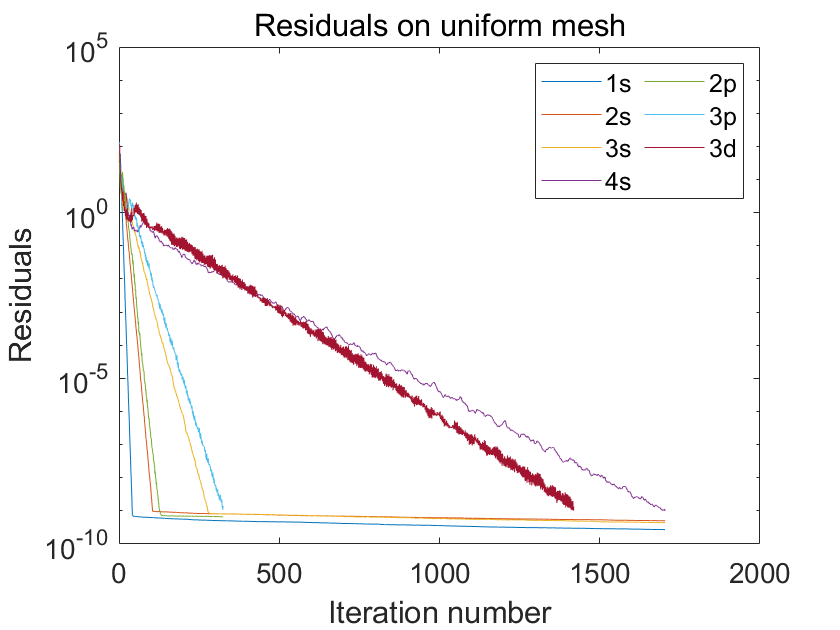}~
    \includegraphics[width=0.45\linewidth]{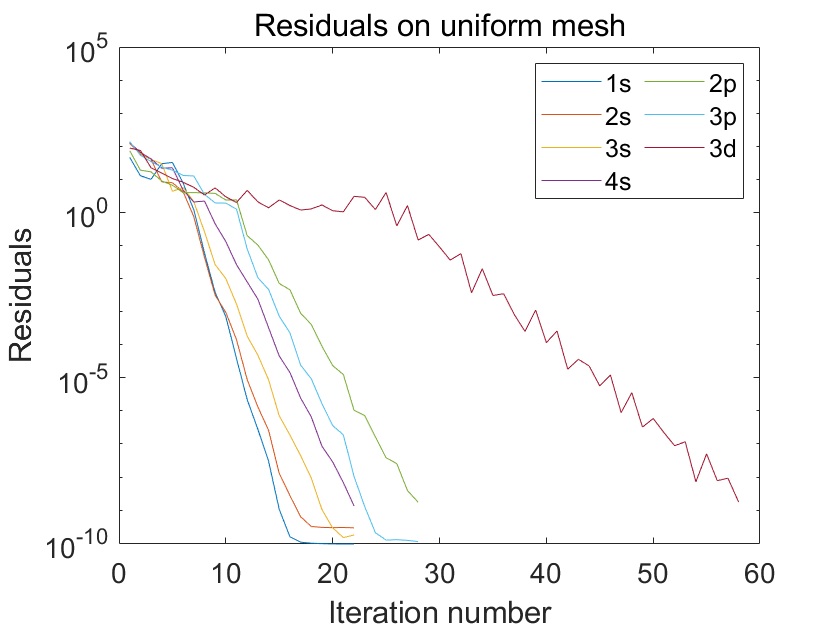}
    \includegraphics[width=0.45\linewidth]{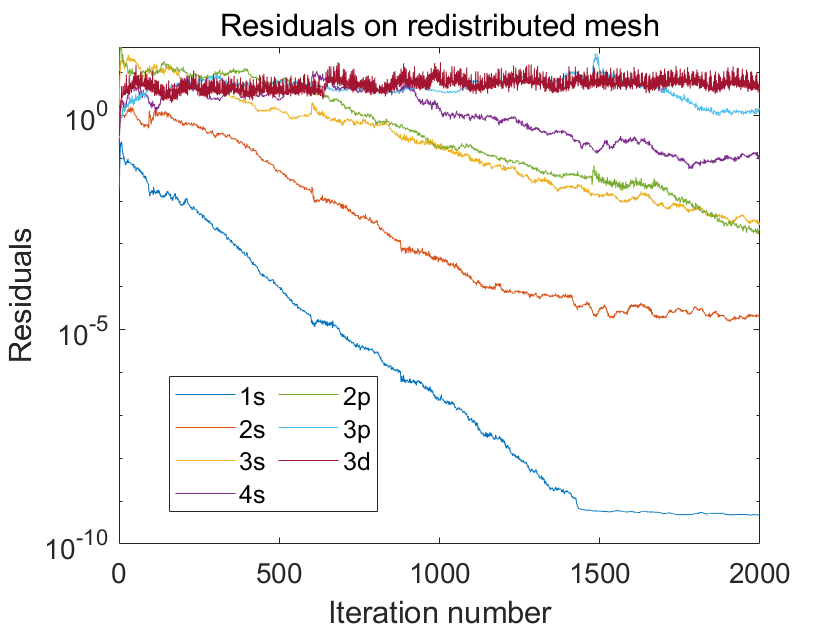}~
    \includegraphics[width=0.45\linewidth]{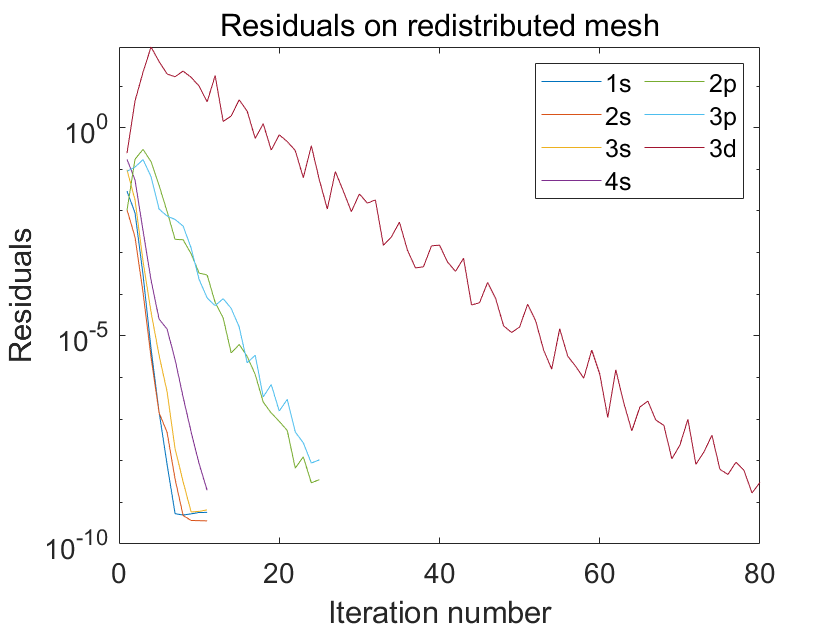}
    \caption{Iteration numbers for LOBPCG in solving the iron atom on a uniform mesh (top) and redistributed mesh (bottom). The left column displays the results without a preconditioner, while the right column shows the results with a preconditioner \eqref{eq:precon}.}
    \label{fig:Fe-precon}
\end{figure}

To validate the effectiveness of the preconditioner, we compare the LOBPCG method with and without it using the iron example. The implementation of LOBPCG is based on the matlab package \texttt{BLOPEX} \cite{knyazev2007block}.  We assess performance on two meshes during the first SCF iteration: a uniform mesh and a redistributed mesh, with a Lagrange polynomial degree of 4 and 80 elements over the interval $[0, 20]$. The redistributed mesh is generated by applying the moving mesh method to the uniform mesh. The stopping criterion for LOBPCG is set when the residual falls below $1.0e-9$. The results are shown in \Cref{fig:Fe-precon}. From the top left of \Cref{fig:Fe-precon}, we observe that on the uniform mesh, more than 1000 iterations are required to achieve convergence without the preconditioner. However, with the preconditioner \eqref{eq:precon}, the iteration number for all orbitals drops below 50. The situation worsens on the redistributed mesh without a preconditioner, as seen in the bottom left of \Cref{fig:Fe-precon}, where only the 1s orbital converges within 2000 steps, while the others do not. This occurs because the mesh sizes on the redistributed mesh vary widely, leading the Hamiltonian matrix to become ill-conditioned. In contrast, after applying the preconditioner, convergence is achieved within 80 steps, significantly speeding up the algorithm.

\subsubsection{Moving mesh method with high-order finite elements}
\begin{figure}[!ht]
    \centering
    \includegraphics[width=0.45\linewidth]{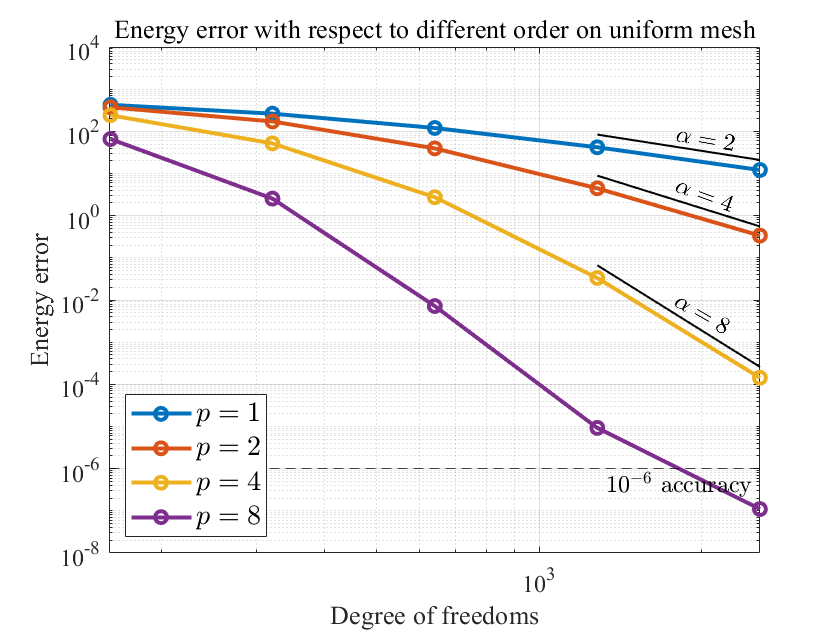}~    
    \includegraphics[width=0.45\linewidth]{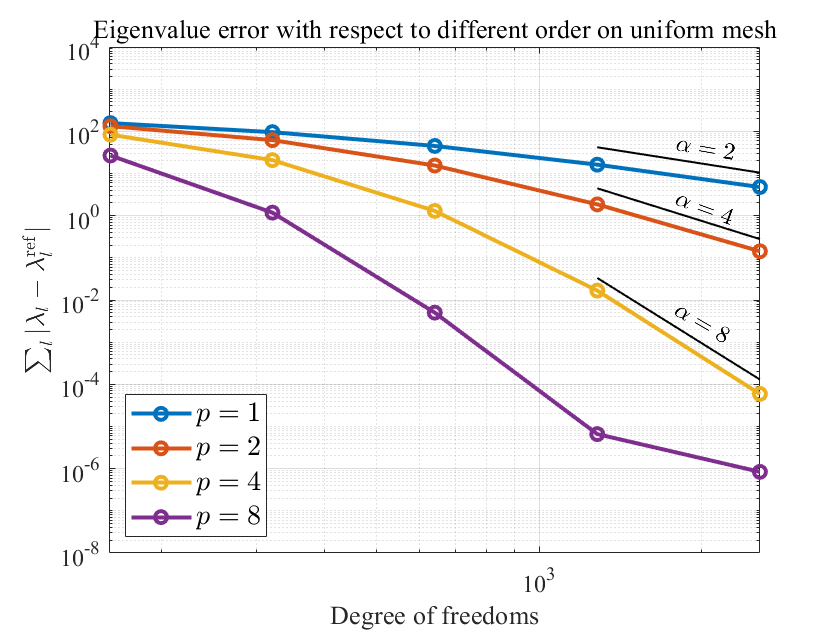}~
    \caption{Convergence results for iron atom on a uniform mesh with respect to different orders. Left is for the total energy and right is for the summation of the eigenvalues.}
    \label{fig:Fe-uniform}
\end{figure}

We first highlight the substantial benefits of high-order elements compared to low-order elements. A series of experiments are conducted using uniform meshes with order $p=1,2,4,8$. The results are displayed in \Cref{fig:Fe-uniform}. The referenced value is from the NIST database. Note that in the database, it only shows the value with 6 digits, hence the NIST result is reproduced when the error is less than $10^{-6}$. The left figure in \Cref{fig:Fe-uniform} displays the energy error with respect to the number of degree of freedoms (DOFs). It is found that the convergence rate for $p=1,2,4$ agrees with the theoretical rate $o(h^{2p})$. The convergence rate for $p=8$ is not plotted as the reference data is accurate up to $10^{-6}$ accuracy. Additionally, on the finest mesh with $2,561$ DOFs, the methods using $p=1,2$ and $4$ elements fail to achieve the desired accuracy. Specifically, for $p=4$, there are $640$ elements on the finest mesh, and the energy error is around $10^{-4}$.
whereas only $p=8$ meets the target accuracy using $320$ elements. Such a phenomenon indicates that the importance of using high-order elements. A similar observation can be drawn from the eigenvalue errors shown in the right of \Cref{fig:Fe-uniform}.

\begin{figure}[!ht]
    \centering
    \includegraphics[width=0.45\linewidth]{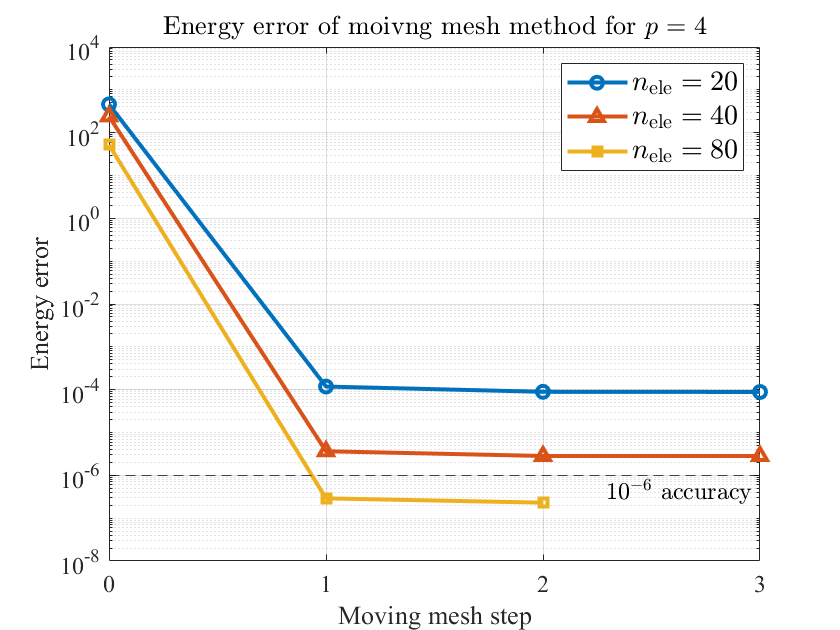}~
    \includegraphics[width=0.45\linewidth]{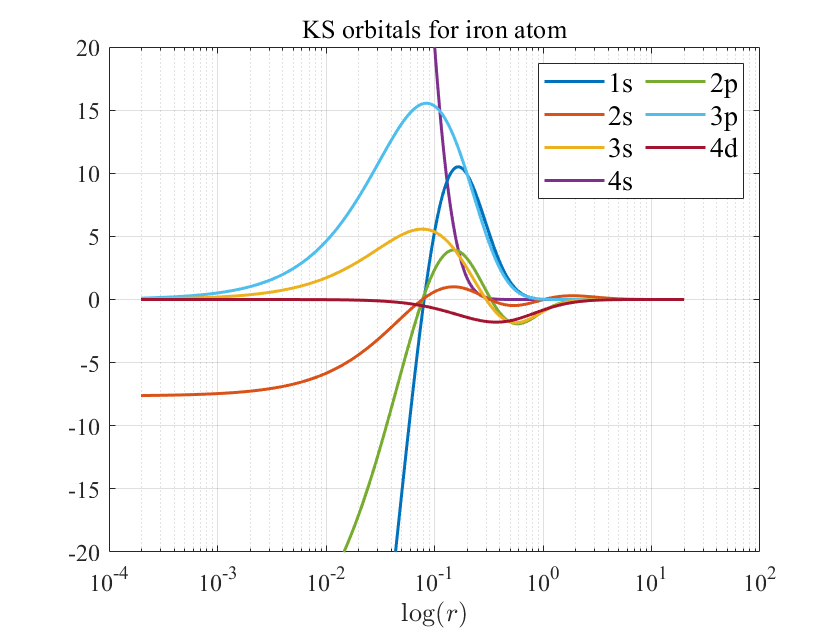}
    \caption{Moving mesh method for the iron atom with order $p=4$.}
    \label{fig:Fe-mmm}
\end{figure}

However, even for the $p=8$ finite element method, over hundreds elements are required to reproduce the NIST result of the iron atom. To reduce the number of required elements, the moving mesh method is applied. The moving mesh method using order $p=4$ finite element for the iron is studied, as shown in \Cref{fig:Fe-mmm}. We test the method with number of elements $n_\mathrm{ele}=20,40$ and $80$. The energy errors with respect to the moving mesh step are displayed in the left of \Cref{fig:Fe-mmm}. The results show that $10^{-6}$ accuracy is achieved when $n_\mathrm{ele}=80$, which attains higher accuracy with far fewer elements than the results obtained from the uniform mesh. Note that the moving mesh method is based on equi-distributing the monitor function \eqref{eq:monitor}. We then plot the KS orbitals as demonstrated in right of \Cref{fig:Fe-mmm}. For a clear demonstration of the variations of the orbitals, the logscale for the $x$-axis is used. It can be found that for the distance to nucleus around $5\times10^{-2}$  to $10^0$, the variations for orbitals are large. Thus, to obtain an accurate result, an adequate number of mesh grids are required. Furthermore, for the region $r>10^{0}$, the orbitals tend to be 0, which implies that few mesh grids are needed in this region. 

\begin{figure}[!ht]
    \centering
    \includegraphics[width=0.45\linewidth]{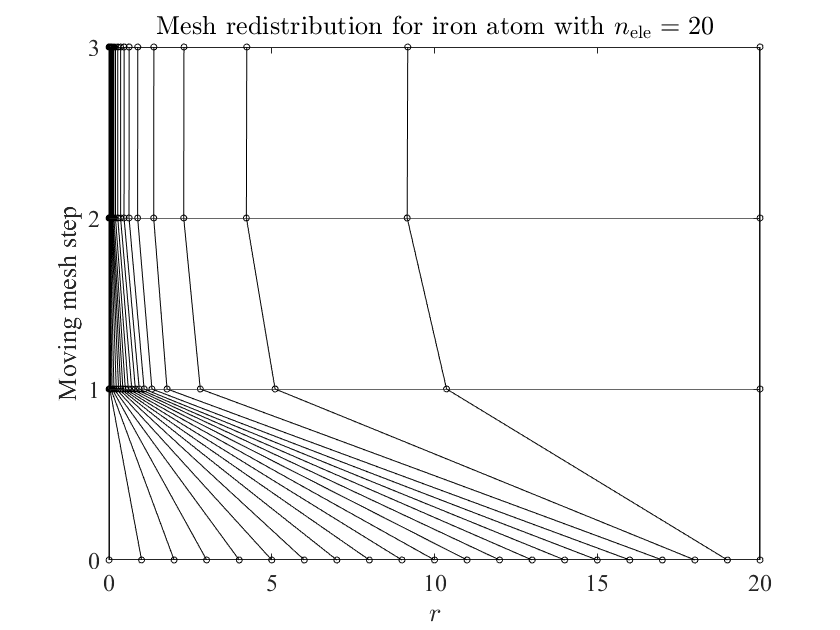}~
    \includegraphics[width=0.45\linewidth]{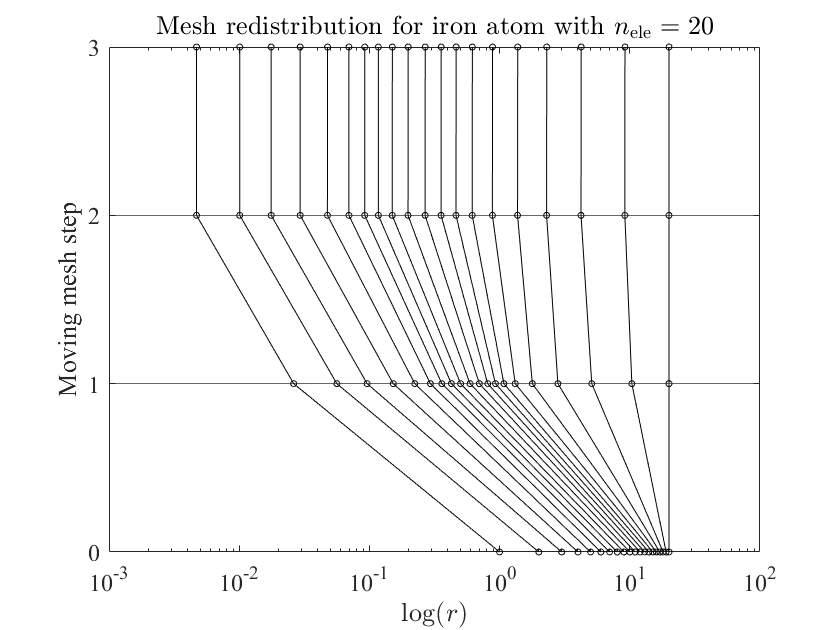}
    \includegraphics[width=0.45\linewidth]{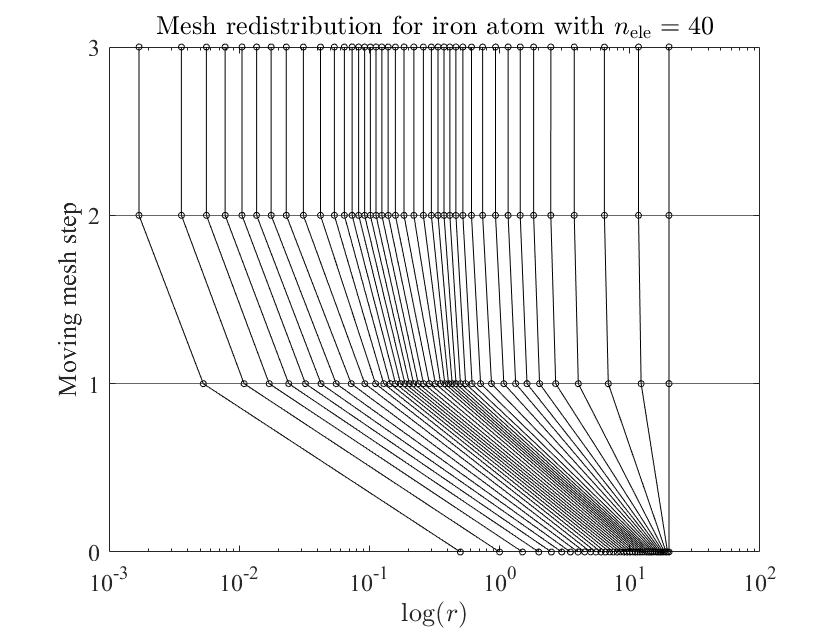}~
    \includegraphics[width=0.45\linewidth]{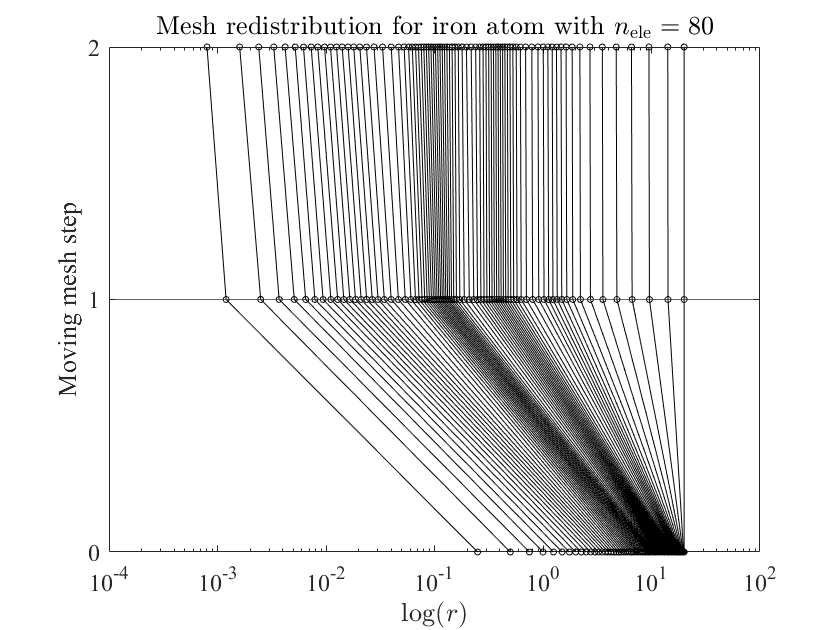}
    \caption{Meshes of iron atom for different number of elements. \label{fig:Fe-mesh}}
    
\end{figure}

The redistribution of the mesh grids using the moving mesh method are displayed in \Cref{fig:Fe-mesh}. The top left of \Cref{fig:Fe-mesh} represent the mesh redistribution for $n_\mathrm{ele}=20$. Obviously, it is found the mesh grids move towards the origin, i.e., the nuclei position. To better observe the distribution around the nucleus, we again plot the meshes using a logarithmic scale on the $x$-axis, as indicated in the top right of \Cref{fig:Fe-mesh}. From this figure, we found that most mesh grids are located at $r<10^0$. This observation is further confirmed for meshes with  $n_\mathrm{ele}=40$ and $80$, as shown in the bottom of \Cref{fig:Fe-mesh}. For all three simulations, a maximum of $3$ moving mesh steps is required to reach convergence, indicating that the redistributed mesh can be obtained with minimal computational cost. It is also noted that for the redistributed mesh using the moving mesh method, the region with the highest concentration of mesh grids is $[5\times10^{-2},10^0]$, particularly in the case of $n_\mathrm{ele}=80$. This phenomenon differs slightly from the exponential mesh, which is essentially uniform on the logarithmic scale of the  $x$-axis,  but it aligns with the fact mentioned in the previous paragraph  the earlier observation that the orbitals exhibit significant variations in the region $[5\times10^{-2},10^{0}]$.  

\begin{figure}[!ht]
    \centering
    \includegraphics[width=0.45\linewidth]{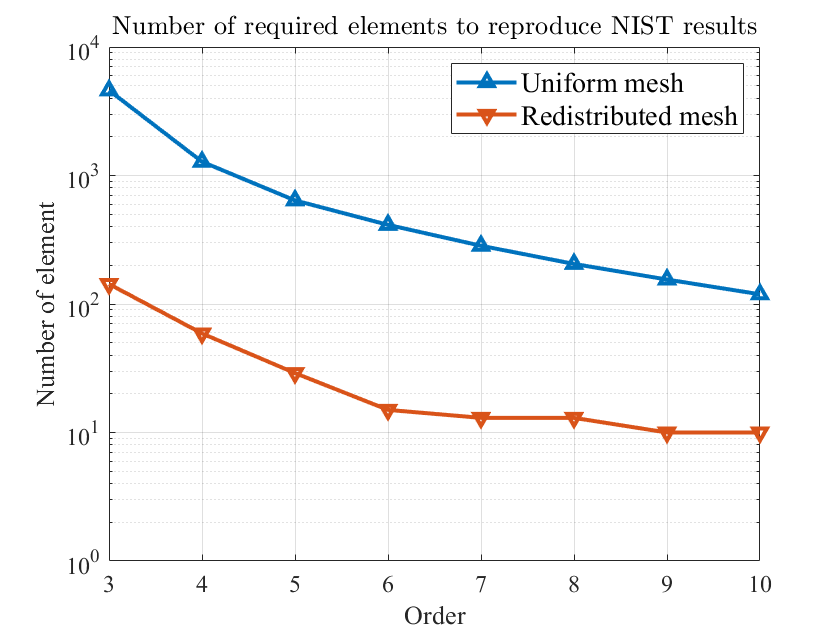}~
    \includegraphics[width=0.45\linewidth]{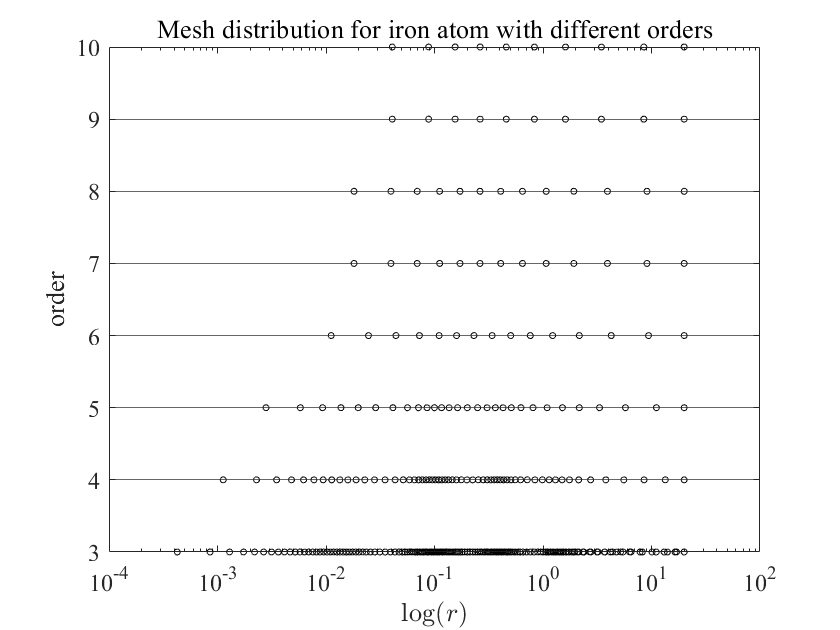}
    \caption{The number of elements of  Fe atom in fixed and moving grids at different orders is compared.\label{fig:Fe-nele}}
\end{figure}

We further show the effectiveness of the high-order elements and the moving mesh method by comparing the number of required elements to match the NIST result. The results and meshes are illustrated in  the left and right of \Cref{fig:Fe-nele}, respectively. It is evident that the number of required elements for the moving mesh method is at most $1/20$ of that required for the uniform mesh when $p\le 7$, highlighting the significant reduction in element count when the moving mesh strategy is utilized. Specifically, only $9$ elements are needed to achieve the same result with NIST database. Additionally, the final redistributed meshes with different order $p$ is plotted in the right of \Cref{fig:Fe-nele}. Similar results on the meshes with the previous discussed case when $p=4$ can be observed for $p\le 5$ . While for  $p>6$, the mesh tends to be like the exponential mesh. 

\subsection{Uranium atom}

To further illustrate the effectiveness and accuracy of the presented method, the uranium atom with atomic number $Z=92$ is studied. In this example, we compare the result with the state-of-art code for atomic structure calculations, i.e., \texttt{featom} \cite{certik2024highorder}, in which a result up to $10^{-9}$ accuracy is given.  

\begin{figure}[!ht]
    \centering
    \includegraphics[width=0.32\linewidth]{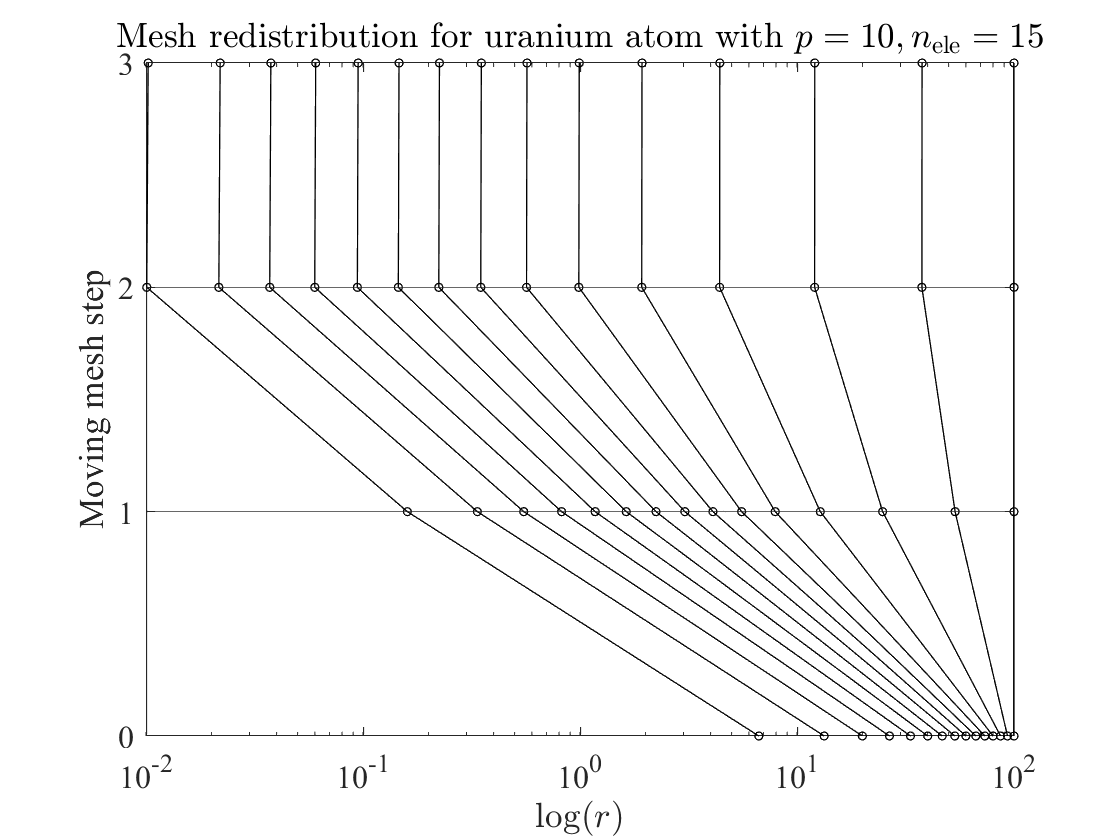}  \includegraphics[width=0.32\linewidth]{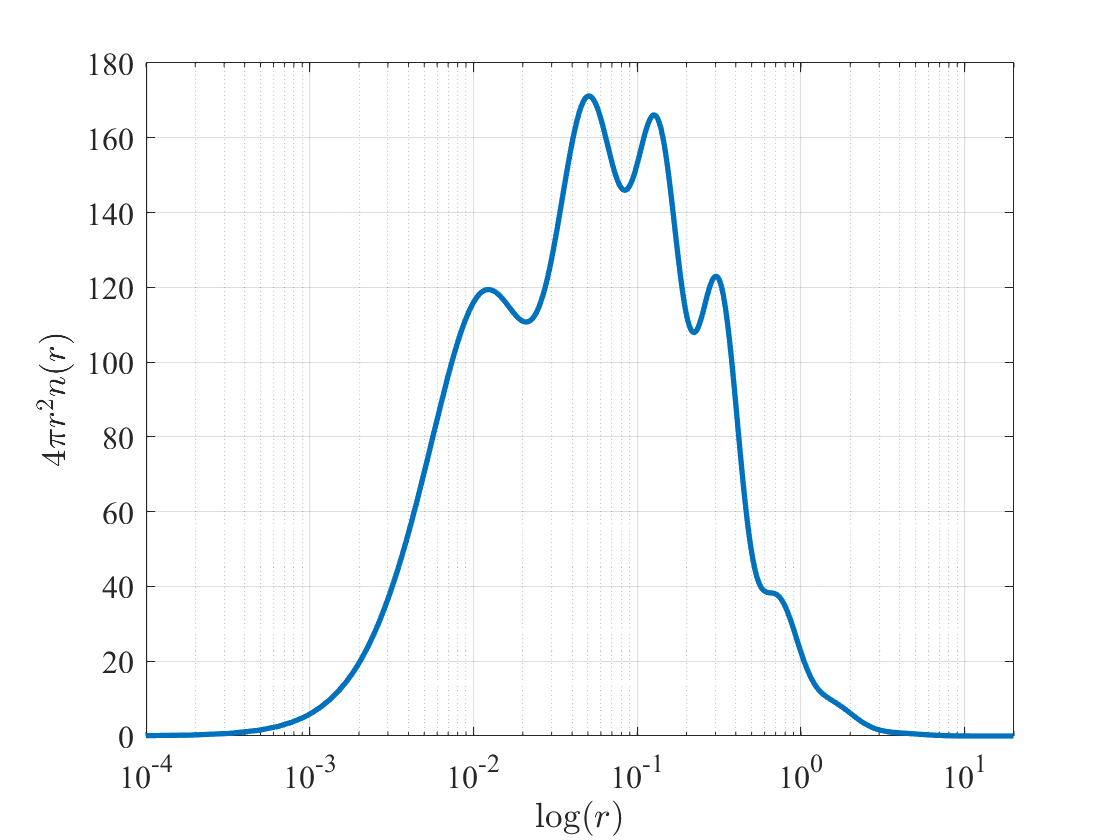}~    \includegraphics[width=0.32\linewidth]{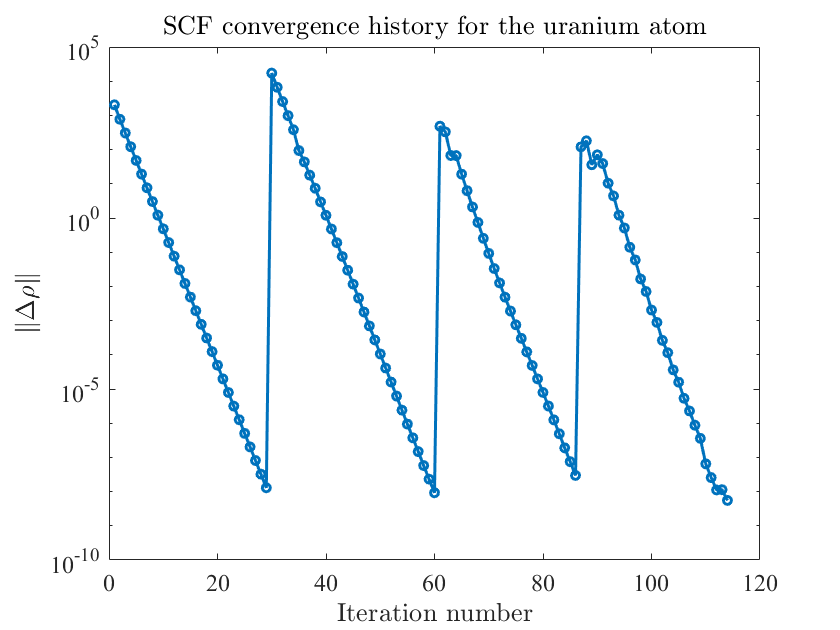}
    \caption{(a) the grid movement distribution map for uranium atom.(b) electron density distribution map for uranium atom .(c) The convergence history for uranium atom. \label{fig:u}}
\end{figure}

Due to the effectiveness of the high-order elements, as mentioned in the previous subsection, we adopt $p=10$ in this example. Furthermore, the computational domain is set as large as $[0,100]$, to avoid the error introduced from the domain size \cite{certik2024highorder}. It is found that a large increment of the domain size will not cause a large number of elements, which can be controlled in a very qualified way by the moving mesh method. In this example, $n_\mathrm{ele}$ is chosen as $15$. The mesh redistribution process is displayed in the left of \Cref{fig:u}. Three moving mesh steps are employed to achieve the convergence of the mesh. It is observed that only two elements are in the right of $10^1$, which implies that the increasing of the domain size will not cause the large increase of the element. The electron density is plotted in the middle of \Cref{fig:u}, in which the peaks are clearly observed in the region $[10^{-2},10^{0}]$. The SCF convergence is displayed in the right of \Cref{fig:u}, in which a linear mixing scheme is adopted.  

The energy and eigenvalue errors are listed in \Cref{tab:u}. It is found that for $n_\mathrm{ele}=15$, we are able to achieve the result with accuracy up to $10^{-9}$ Hartree. The errors of the eigenvalues also reach such an accuracy.
\begin{table}[!ht]
    \centering
    \caption{Eigenvalues for the uranium atom.}\label{tab:u}
    \begin{tabular}{r|r|r|c}     
    \toprule  uranium&$E_\mathrm{tot}$&$E_\mathrm{tot}^{\mathrm{Ref}}$&Error\\
    \midrule
    -& -25658.41788885&	-25658.41788885	&3.36E-09	 \\ \midrule
    $nl$	& $\varepsilon_{nl}$	&$\varepsilon_n^{\mathrm{Ref}}$	&Error\\ \midrule
1	&-3689.35513983&	-3689.35513984&	1.02E-08\\
2&	-639.77872808&	-639.77872809&	7.33E-09\\
3&	-619.10855018&	-619.10855018&	2.53E-09\\
4&	-161.11807321&	-161.11807321&	4.46E-09\\
5&	-150.97898016&	-150.97898016&	6.90E-10\\
6&	-131.97735828&	-131.97735828&	1.14E-09\\
7&	-40.52808424&	-40.52808425&	9.36E-09\\
8&	-35.85332083&	-35.85332083&	2.00E-09\\
9&	-27.12321229&	-27.12321230&	5.03E-09\\
10&	-15.02746006&	-15.02746007&	5.54E-09\\
11&	-8.82408940&	-8.82408940&	3.52E-09\\
12&	-7.01809220&	-7.01809220&	5.14E-10\\
13&	-3.86617513&	-3.86617513&	5.75E-11\\
14&	-1.32597631&	-1.32597632&	6.58E-09\\
15&	-0.82253797&	-0.82253797&	3.96E-09\\
16&	-0.36654335&	-0.36654335&	1.63E-09\\
17&	-0.14319018&	-0.14319018&	1.86E-09\\
18&	-0.13094786&	-0.13094786&	3.90E-10\\
\bottomrule
    \end{tabular}

\end{table}

\subsection{Ground state energy}

\begin{figure}[!h]
    \centering
    \scalebox{0.3}{
    \input{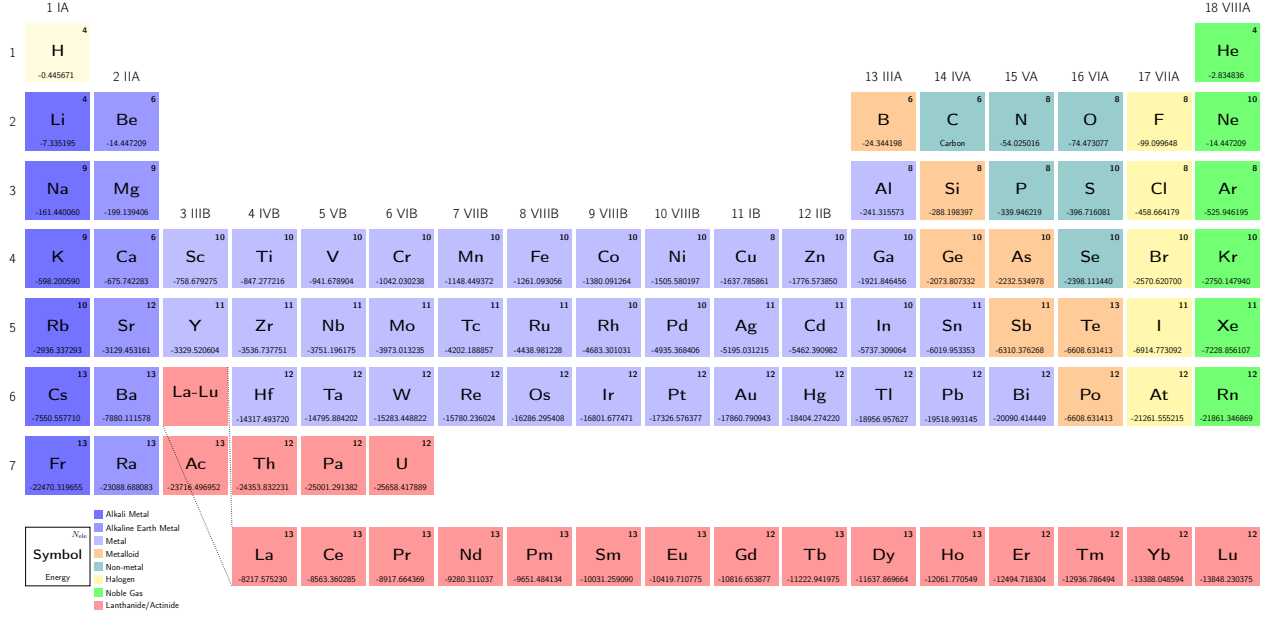}
    } 
    \caption{The number of required $p=10$ finite elements to reproduce the results from NIST database \cite{kotochigova2009atomic}. The atomic number is ranging from $1$ to $92$. The computational domain  is set as $[0,20]$  for atoms in the first four rows and $[0,100]$ for the remained atoms.     \label{fig:periodic-table}}
\end{figure}

\begin{figure}[!h]
    \centering
    \includegraphics[width=0.32\linewidth]{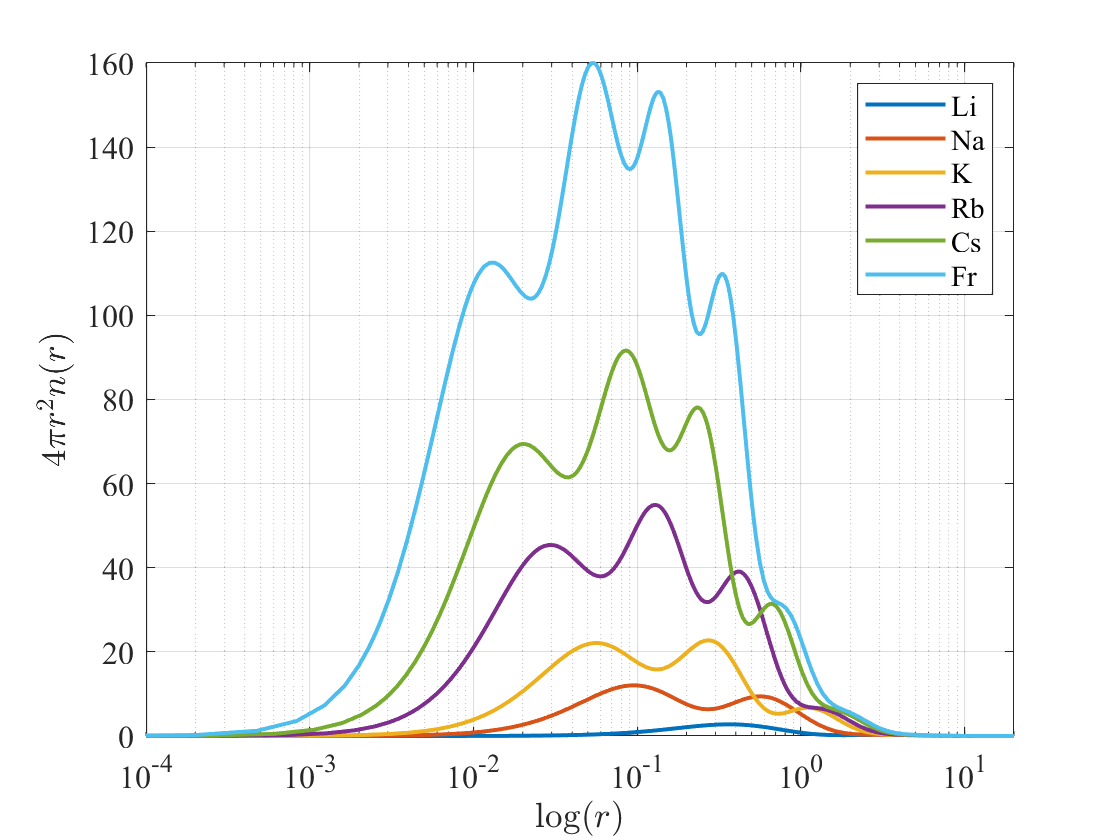}~        \includegraphics[width=0.32\linewidth]{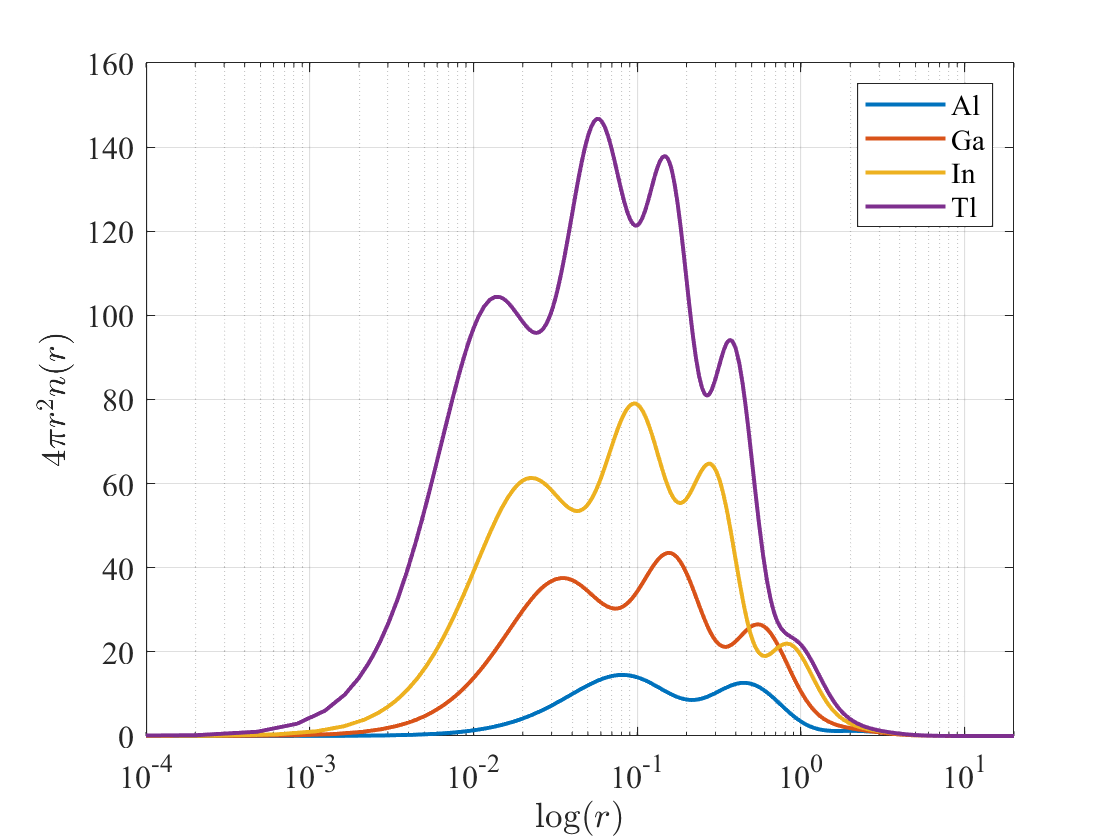}~
    \includegraphics[width=0.32\linewidth]{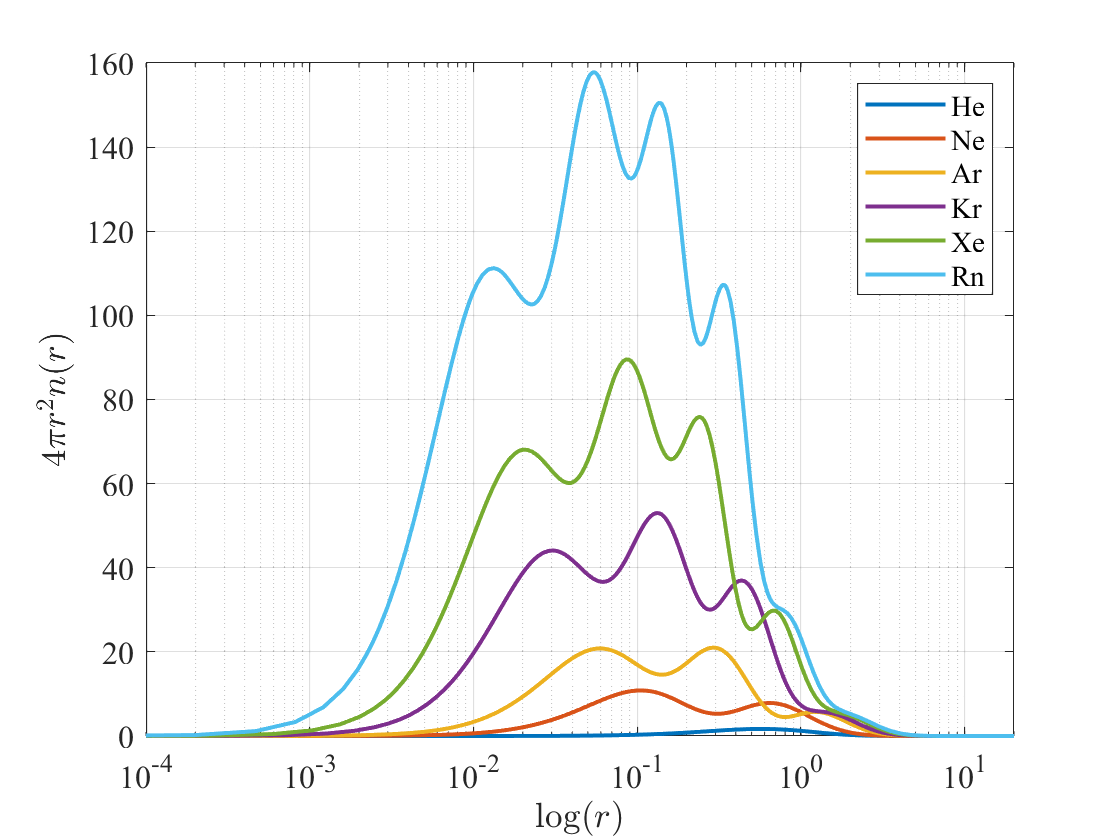}
    \caption{Electron densities for atoms in Group IA (left), IIIA (middle), and VIIIA (right).}
    \label{fig:allrho}
\end{figure}

Finally, we list the minimal number of required number of elements with order $p=10$ using the presented method to reproduce the results from NIST database \cite{kotochigova2009atomic} for all the atoms with atomic numbers from $1$ to $92$. The results are listed in \Cref{fig:periodic-table}. The computational domain for the first four rows in the periodic table is set as $[-20,20]$, and it is set as $[0,100]$ for the atoms in the last three rows. As a result, at most $10$ finite elements are required to achieve the desired accuracy for the first four rows, and at most $13$ finite elements are  required for the remaining atoms. Such results show the accuracy and efficiency of the presented method. The electron densities for atoms in Group IA, IIIA, and VIIIA are displayed in \Cref{fig:allrho}. 

\section{Conclusion}
This paper presents a high-order accurate moving mesh method for solving the radial Kohn--Sham equation. Compared to a uniform mesh, the moving mesh approach successfully reduces the number of elements to at most $1/20$ for the iron atom. For the uranium atom, only $15$ elements are needed to achieve an accuracy of $10^{-9}$ Hartree. Notably, the method is parameter-free and applicable to all atoms under the same numerical configuration, demonstrating both its generality and robustness. Additionally, only a few mesh adaptation steps are required to achieve convergence, ensuring that the computational cost remains low. Using this method, we reproduce the NIST database results for atoms with atomic numbers ranging from 1 to 92 within just 13 elements. The proposed approach offers an efficient and parameter-free solution for the radial Kohn--Sham equation, making it well-suited for generating pseudopotentials or linearized muffin-tin orbitals, which we plan to explore in future work.

\section*{Acknowledgement}
This work is supported by the National Natural Science Foundation of China (Nos. 12201130 and 12326362) and the Guangzhou Municipal Science and Technology Bureau (No. 2023A04J1321).

\bibliographystyle{plain}
\bibliography{main.bib}

\begin{thebibliography}{10}

\bibitem{bao2012hadaptive}
Gang Bao, Guanghui Hu, and Di~Liu.
\newblock An $h$-adaptive finite element solver for the calculations of the
  electronic structures.
\newblock {\em Journal of Computational Physics}, 231(14):4967--4979, May 2012.

\bibitem{castro2006octopus}
Alberto Castro, Heiko Appel, Micael Oliveira, Carlo~A Rozzi, Xavier Andrade,
  Florian Lorenzen, Miguel~AL Marques, EKU Gross, and Angel Rubio.
\newblock Octopus: a tool for the application of time-dependent density
  functional theory.
\newblock {\em Physica Status Solidi (b)}, 243(11):2465--2488, 2006.

\bibitem{certik2024highorder}
Ond{\v r}ej {\v C}ert{\'i}k, John~E. Pask, Isuru Fernando, Rohit Goswami,
  N.~Sukumar, {\relax Lee}.~A. Collins, Gianmarco Manzini, and Ji{\v r}{\'i}
  Vack{\'a}{\v r}.
\newblock High-order finite element method for atomic structure calculations.
\newblock {\em Computer Physics Communications}, 297:109051, April 2024.

\bibitem{certik2013dftatom}
Ond{\v r}ej {\v C}ert{\'i}k, John~E. Pask, and Ji{\v r}{\'i} Vack{\'a}{\v r}.
\newblock Dftatom: {{A}} robust and general {{Schr{\"o}dinger}} and {{Dirac}}
  solver for atomic structure calculations.
\newblock {\em Computer Physics Communications}, 184(7):1777--1791, July 2013.

\bibitem{de1973good}
Carl de~Boor.
\newblock Good approximation by splines with variable knots.
\newblock In {\em Spline Functions and Approximation Theory: Proceedings of the
  Symposium held at the University of Alberta, Edmonton May 29 to June 1,
  1972}, pages 57--72. Springer, 1973.

\bibitem{fiolhais2008primer}
Carlos Fiolhais, Fernando Nogueira, and Miguel~AL Marques.
\newblock {\em A primer in density functional theory}, volume 620.
\newblock Springer, 2008.

\bibitem{gonze2009abinit}
Xavier Gonze, Bernard Amadon, P-M Anglade, J-M Beuken, Fran{\c{c}}ois Bottin,
  Paul Boulanger, Fabien Bruneval, Damien Caliste, Razvan Caracas, Michel
  C{\^o}t{\'e}, et~al.
\newblock {ABINIT: First-principles approach to material and nanosystem
  properties}.
\newblock {\em Computer Physics Communications}, 180(12):2582--2615, 2009.

\bibitem{gulans2014exciting}
Andris Gulans, Stefan Kontur, Christian Meisenbichler, Dmitrii Nabok, Pasquale
  Pavone, Santiago Rigamonti, Stephan Sagmeister, Ute Werner, and Claudia
  Draxl.
\newblock Exciting: A full-potential all-electron package implementing
  density-functional theory and many-body perturbation theory.
\newblock {\em Journal of Physics: Condensed Matter}, 26(36):363202, September
  2014.

\bibitem{hohenberg1964inhomogeneous}
Pierre Hohenberg and Walter Kohn.
\newblock Inhomogeneous electron gas.
\newblock {\em Physical Review}, 136(3B):B864, 1964.

\bibitem{huang2011adaptive}
Weizhang Huang and Robert~D. Russell.
\newblock {\em Adaptive {{Moving Mesh Methods}}}, volume 174 of {\em Applied
  {{Mathematical Sciences}}}.
\newblock Springer New York, New York, NY, 2011.

\bibitem{Adaptive-finite-element3}
Andrew~V Knyazev.
\newblock Toward the optimal preconditioned eigensolver: Locally optimal block
  preconditioned conjugate gradient method.
\newblock {\em SIAM Journal on Scientific Computing}, 23(2):517--541, 2001.

\bibitem{knyazev2007block}
Andrew~V Knyazev, Merico~E Argentati, Ilya Lashuk, and Evgueni~E Ovtchinnikov.
\newblock {Block locally optimal preconditioned eigenvalue xolvers (BLOPEX) in
  HYPRE and PETSc}.
\newblock {\em SIAM Journal on Scientific Computing}, 29(5):2224--2239, 2007.

\bibitem{kohn1965self}
Walter Kohn and Lu~Jeu Sham.
\newblock Self-consistent equations including exchange and correlation effects.
\newblock {\em Physical Review}, 140(4A):A1133, 1965.

\bibitem{kotochigova1997local}
Svetlana Kotochigova, Zachary~H Levine, Eric~L Shirley, Mark~D Stiles, and
  Charles~W Clark.
\newblock {Local-density-functional calculations of the energy of atoms}.
\newblock {\em Physical Review A}, 55(1):191, 1997.

\bibitem{kotochigova2009atomic}
Svetlana Kotochigova, Zachary~H Levine, Eric~L Shirley, Mark~D Stiles, and
  Charles~W Clark.
\newblock Atomic reference data for electronic structure calculations (version
  1.6), 2009.
\newblock Available at
  \url{https://www.nist.gov/pml/atomic-reference-data-electronic-structure-calculations/atomic-reference-data-electronic-7},
  Accessed: 2024-08-01.

\bibitem{lehtola2019fully}
Susi Lehtola.
\newblock Fully numerical {{Hartree}}-{{Fock}} and density functional
  calculations. {{I}}. {{Atoms}}.
\newblock {\em International Journal of Quantum Chemistry}, 119(19):e25945,
  October 2019.

\bibitem{lehtola2020fully}
Susi Lehtola.
\newblock Fully numerical calculations on atoms with fractional occupations and
  range-separated exchange functionals.
\newblock {\em Physical Review A}, 101(1):012516, January 2020.

\bibitem{oliveira2008generating}
Micael~J.T. Oliveira and Fernando Nogueira.
\newblock Generating relativistic pseudo-potentials with explicit incorporation
  of semi-core states using {{APE}}, the {{Atomic Pseudo-potentials Engine}}.
\newblock {\em Computer Physics Communications}, 178(7):524--534, April 2008.

\bibitem{oulne2011variation}
M~Oulne.
\newblock Variation and series approach to the thomas--fermi equation.
\newblock {\em Applied Mathematics and Computation}, 218(2):303--307, 2011.

\bibitem{pickett1989pseudopotential}
Warren~E Pickett.
\newblock Pseudopotential methods in condensed matter applications.
\newblock {\em Computer Physics Reports}, 9(3):115--197, 1989.

\bibitem{romanowski2007numerical}
Zbigniew Romanowski.
\newblock Numerical solution of {K}ohn-{S}ham equation for atom.
\newblock {\em Acta Physica Polonica B}, 38(10), 2007.

\bibitem{romanowski2008bspline}
Zbigniew Romanowski.
\newblock A {{B-spline}} finite element solution of the {{Kohn}}--{{Sham}}
  equation for an atom.
\newblock {\em Modelling and Simulation in Materials Science and Engineering},
  16(1):015003, January 2008.

\bibitem{romanowski2011bspline}
Zbigniew Romanowski.
\newblock B-spline solver for one-electron {{Schr{\"o}dinger}} equation.
\newblock {\em Molecular Physics}, 109(22):2679--2691, November 2011.

\bibitem{singh2006planewaves}
David~J Singh and Lars Nordstrom.
\newblock {\em {Planewaves, Pseudopotentials, and the LAPW method}}.
\newblock Springer Science \& Business Media, 2006.

\bibitem{skriver2012lmto}
Hans~L Skriver.
\newblock {\em {The LMTO method: muffin-tin orbitals and electronic
  structure}}, volume~41.
\newblock Springer Science \& Business Media, 2012.

\bibitem{tang2005moving}
Tao Tang.
\newblock Moving mesh methods for computational fluid dynamics.
\newblock In Z.-C. Shi, Z.~Chen, T.~Tang, and D.~Yu, editors, {\em Contemporary
  {{Mathematics}}}, volume 383, pages 141--173. American Mathematical Society,
  Providence, Rhode Island, 2005.

\bibitem{uzulis2022radial}
J{\=a}nis U{\v z}ulis and Andris Gulans.
\newblock Radial {{Kohn}}--{{Sham}} problem via integral-equation approach.
\newblock {\em Journal of Physics Communications}, 6(8):085002, August 2022.

\bibitem{vosko1980accurate}
Seymour~H Vosko, Leslie Wilk, and Marwan Nusair.
\newblock Accurate spin-dependent electron liquid correlation energies for
  local spin density calculations: a critical analysis.
\newblock {\em Canadian Journal of Physics}, 58(8):1200--1211, 1980.

\end{thebibliography}
\end{document}